%% file: adap_res.tex
\documentclass[11pt]{article}
\usepackage{fullpage,graphicx,psfrag,color,verbatim,url,algorithm,algorithmic}
\input defs.tex

\title{Adaptive Restart for Accelerated Gradient Schemes}
\author{Brendan O'Donoghue \and Emmanuel Cand\`{e}s}
\begin{document} \maketitle

\begin{abstract}
In this paper we demonstrate a simple heuristic adaptive restart technique that can
dramatically improve the convergence rate of accelerated gradient schemes.  The
analysis of the technique relies on the observation that these schemes exhibit
two modes of behavior depending on how much momentum is applied. In what we
refer to as the `high momentum' regime the iterates generated by an accelerated
gradient scheme exhibit a periodic behavior, where the period is
proportional to the square root of the local condition number of the objective
function.  This suggests a restart technique whereby we reset the momentum
whenever
we observe periodic behavior. We provide analysis to show that in many cases
adaptively restarting allows us to recover the optimal rate of convergence
with no prior knowledge of function parameters.
\end{abstract}
\section{Introduction} 
%


Accelerated gradient schemes were first
proposed by Yurii Nesterov in 1983, \cite{Nes:83}. He demonstrated a simple
modification to gradient descent that could obtain provably optimal performance
for the complexity class of first-order algorithms applied to minimize smooth
convex functions.  The method, and its successors, are often referred to as
`accelerated methods'. In recent years
there has been a resurgence of interest in first-order optimization methods
\cite{tfocs, nes:07,tseng:08,auslender:06, lan:09}, driven primarily by the
need to solve very large problem instances unsuited to second-order methods.

Accelerated gradient schemes can be thought of as \emph{momentum} methods, in
that the step taken at the current iteration depends on the previous
iterations, and where the momentum grows from one iteration to the next. When
we refer to \emph{restarting} the algorithm we mean starting the algorithm
again, taking the current iteration as the new starting point. This erases the
memory of previous iterations and resets the momentum back to zero.

Unlike gradient descent, accelerated methods are not guaranteed to be monotone
in the objective value. A common observation when running an accelerated method
is the appearance of ripples or bumps in the trace of the objective value; these
are seemingly regular increases in the objective, see Figure
(\ref{f-qeff}) for an example. In this paper we demonstrate that this
behavior occurs when the momentum has exceeded a critical value (the optimal
momentum value derived by Nesterov in \cite{nesbook}) and that the period of these
ripples is proportional to the square-root of the (local) condition number of
the function. Separately, we show that the optimal restart interval is also
proportional to the square root of the condition number.  Combining these
results we show that restarting when we observe an increase in the function
value allows us to recover the optimal linear convergence rate in many cases.
Indeed if the function is locally well-conditioned we can use restarting to
obtain a linear convergence rate inside the well-conditioned region.

\paragraph{Smooth unconstrained optimization.} 
We wish to minimize a smooth convex function of a variable
$x \in \reals^n$ \cite{Boyd04}, 
\begin{equation} 
\label{e-optprob}
\begin{array}{ll} 
\mbox{minimize} & f(x) 
\end{array} 
\end{equation}
where $f:\reals^n \rightarrow \reals$ has a Lipschitz continuous
gradient with constant $L$, \ie,
\[ 
\|\nabla f(x) -\nabla f(y)\|_2 \leq L \|x-y\|_2, \quad \forall x,y \in \reals^n.
\]
We shall denote by $f^\star$ the optimal value of the above optimization
problem, if the minimizer exists and is unique then we shall write it as $x^\star$.
Further, a function is said to be strongly convex if
there exists a $\mu > 0$ such that 
\[ 
f(x) \geq f^\star
+(\mu/2)\|x-x^{\star}\|_2^2, \quad \forall x \in \reals^n,
\]
where $\mu$ is referred to as the strong convexity parameter.
The \emph{condition number} of a smooth, strongly convex function is $L/\mu$.

\section{Accelerated methods} 
Accelerated first-order methods to solve (\ref{e-optprob}) were first developed
by Nesterov \cite{Nes:83}, this scheme is from \cite{nesbook}: 
\begin{algorithm}[H]
\caption{Accelerated scheme I}
\label{a-scheme1}
\begin{algorithmic}[1]
\REQUIRE $x^{0} \in \reals^n$, $y^{0} = x^{0}$, $\theta_0=1$ and $q \in [0,1]$
\FOR {$k=0,1,\ldots $}
\STATE $x^{k+1} = y^{k} - t_k \nabla f(y^{k})$ 
\STATE $\theta_{k+1}$ solves $\theta_{k+1}^2 =
(1-\theta_{k+1})\theta_k^2+q\theta_{k+1}$ 
\STATE $\beta_{k+1} = \theta_k(1 -
\theta_k)/(\theta_k^2 + \theta_{k+1})$ 
\STATE $y^{k+1} = x^{k+1} + \beta_{k+1}(x^{k+1} - x^{k})$
\ENDFOR
\end{algorithmic}
\end{algorithm} 
There are many variants of the
above scheme, see, $\eg$, \cite{tseng:08,nes:07,lan:09,auslender:06,fista}.
Note that by setting $q=1$ in the above scheme
we recover gradient descent.
For a smooth convex function the above scheme converges for any 
$t_k \leq 1/L$; setting $t_k = 1/L$ and $q=0$ obtains a guaranteed
convergence rate of
\begin{equation}
\label{e-ksq}
f(x^{k}) - f^{\star} \leq  \frac{4 L \|x^0 - x^{\star}\|^2}{(k+2)^2}.
\end{equation}
If the function is also strongly convex
with parameter $\mu$, then a choice of $q=\mu/L$ (the reciprocal of the
condition number) will achieve
\begin{equation}
\label{e-lincvg}
f(x^{k}) - f^{\star} \leq L \left(1-\sqrt{\frac{\mu}{L}}\right)^k \|x^0 - x^{\star}\|^2.  
\end{equation}
This is often referred to as \emph{linear convergence}. 
With this convergence rate we can achieve an accuracy of $\epsilon$ in
\begin{equation}
\label{e-ordcvg}
\mathcal{O}\left(\sqrt{\frac{L}{\mu}}\log \frac{1}{\epsilon} \right)
\end{equation}
iterations.

In the case of a strongly convex function the following simpler
scheme obtains the same guaranteed rate of convergence \cite{nesbook}:
\begin{algorithm}[H]
\caption{Accelerated scheme II}
\label{a-scheme2}
\begin{algorithmic}[1]
\REQUIRE $x^{0} \in \reals^n$, $y^{0} = x^{0}$
\FOR {$k=0,1,\ldots $}
\STATE $x^{k+1} = y^{k} - (1/L) \nabla f(y^{k})$
\STATE $y^{k+1} = x^{k+1} + \beta^{\star}(x^{k+1} -x^{k})$
\ENDFOR
\end{algorithmic}
\end{algorithm} 
Where we set 
\begin{equation}
\label{e-bstar}
\beta^{\star}=\frac{1- \sqrt{\mu/L}}{1+\sqrt{\mu/L}}.
\end{equation}
Note that in Algorithm \ref{a-scheme1}, using the optimal choice $q = \mu/L$,
we have that $\beta_k \uparrow \beta^\star$. Taking $\beta_k$ to be a momentum
parameter, then for a strongly convex function $\beta^\star$ is the maximum
amount of momentum we should apply; when we have a value of $\beta$ higher than
$\beta^\star$ we refer to it as `high momentum'. We shall return to this point
later.

The convergence of these schemes is optimal in the sense of the lower
complexity bounds derived by Nemirovski and Yudin in \cite{nemirovski:83}.
However, this convergence is only guaranteed when the function parameters $\mu$
and $L$ are known in advance.

\subsection{Robustness} 
A natural question to ask is how robust are accelerated methods to errors in the
estimates of the Lipschitz constant $L$ and strong convexity parameter
$\mu$? For the case of an unknown Lipschitz constant we can estimate 
the optimal step-size by the use of backtracking; see, \eg, \cite{tseng:08,tfocs}.
Estimating the strong convexity parameter is much more challenging.

\begin{figure}[h]
\begin{center} 
\psfrag{k}[t][t]{\small $k$}
\psfrag{fk-fstar}[b][b]{\small $f(x^k) - f^\star$}
\psfrag{q=0}{\small $q=0$}
\psfrag{q=qs/10}{\small $q=q^\star/10$}
\psfrag{q=qs/3}{\small $q=q^\star/3$}
\psfrag{q=qs}{\small $q=q^\star$}
\psfrag{q=qs*3}{\small $q=3q^\star$}
\psfrag{q=qs*10}{\small $q=10q^\star$}
\psfrag{q=1}{\small $q=1$}
\includegraphics[width=0.7\linewidth]{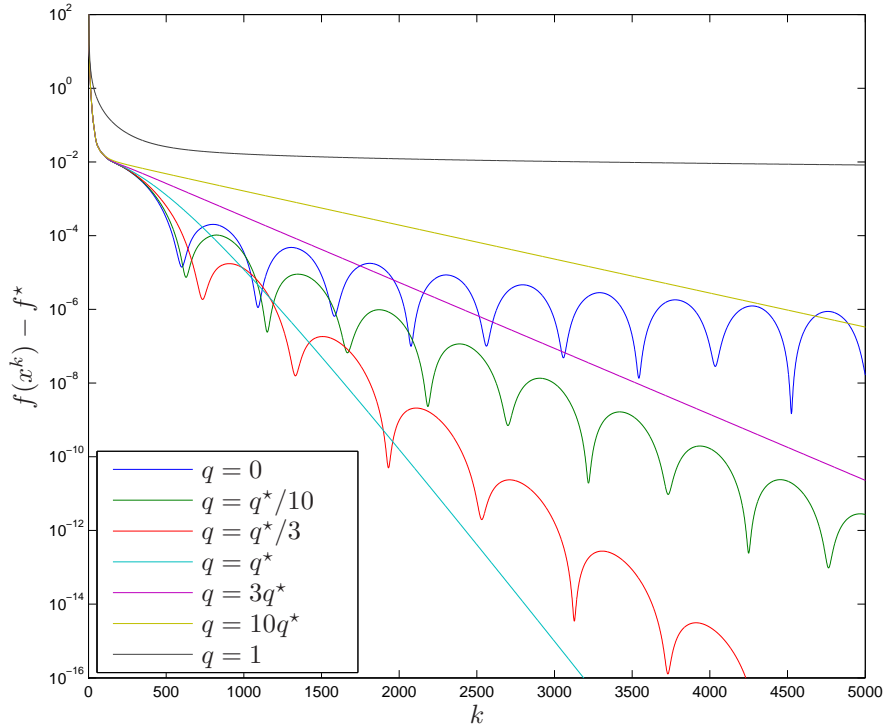} 
\caption{Convergence of Algorithm \ref{a-scheme1} with different estimates of $q$.}
\label{f-qeff} 
\end{center} 
\end{figure}

\paragraph{Estimating the strong convexity parameter.} 
In \cite{nes:07} Nesterov demonstrated a method to bound $\mu$, similar to the
backtracking scheme for $L$ described above. His scheme achieves a
convergence rate quite a bit slower than Algorithm \ref{a-scheme1} with a
known value of $\mu$.
In practice, we often assume or guess that $\mu$ is zero, which corresponds to
setting $q=0$ in Algorithm \ref{a-scheme1}. Indeed many discussions of
accelerated algorithms do not even include a $q$ term; the original
algorithm in \cite{Nes:83} did not use a $q$.
However, this can dramatically slow down the convergence of the algorithm.
Figure \ref{f-qeff} shows Algorithm 
\ref{a-scheme1} applied to minimize a positive definite quadratic function in
$n=200$ dimensions, with optimal choice of $q$ being $q^\star=\mu/L=4.1\times
10^{-5}$ (a condition number of about $2.4\times 10^4$), and step size $t=1/L$.
Each trace is the progress of the algorithm with a different choice of $q$
(hence a different estimate of $\mu$).

We observe that slightly over or underestimating the optimal value of $q$ for the
function can have a severe detrimental effect on the rate of convergence of the
algorithm. We also note the clear difference in behavior between the cases
where we underestimate and where we overestimate $q^\star$; in the latter we observe
monotonic convergence but in the former we notice the appearance of regular
ripples or bumps in the traces.


\begin{figure} 
\begin{center} 
\psfrag{x1}[t][]{$x_1$}
\psfrag{x2}[b][b]{$x_2$}
\psfrag{q = q-opt}{\footnotesize $q=q^\star$}
\psfrag{q = 0}{\footnotesize $q=0$}
\includegraphics[width=0.5\linewidth]{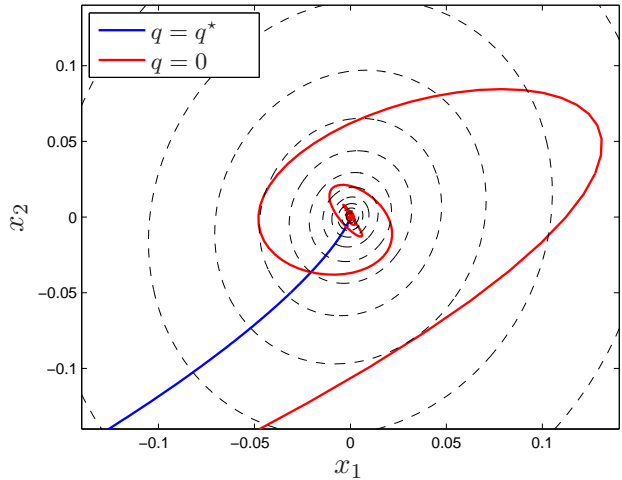} 
\caption{Sequence trajectories under Algorithm \ref{a-scheme1}.}
\label{f-traj_noadap} 
\end{center} 
\end{figure}

\paragraph{Interpretation.} 
The optimal momentum depends on the condition number of the
function; specifically, higher momentum is required when the function has a
higher condition number.  Underestimating the amount of momentum required leads
to slower convergence.  However we are more often in the other regime, that of
overestimated momentum, because generally $q=0$, in which case $\beta_k
\uparrow 1$; this corresponds to high momentum and rippling behavior, as we
see in Figure \ref{f-qeff}. This can be visually understood in Figure
(\ref{f-traj_noadap}), which shows the trajectories of sequences
generated by Algorithm \ref{a-scheme1} minimizing a positive definite quadratic
in two dimensions, under $q = q^\star$, the optimal choice of $q$, and $q=0$.
The high momentum causes the trajectory to overshoot the minimum and oscillate
around it. This causes a rippling in the function values along the trajectory.
Later we shall demonstrate that the period of these ripples is
proportional to the square root of the (local) condition number of the
function.

Lastly we mention that the condition number is a \emph{global} parameter; the
sequence generated by an accelerated scheme may enter regions that are locally
better conditioned, say, near the optimum.  In these cases the choice of $q =
q^\star$ is appropriate outside of this region, but once we enter it we expect
the rippling behavior associated with high momentum to emerge, despite the optimal
choice of $q$. 

\section{Restarting}
\subsection{Fixed restart}
For strongly convex functions an alternative to choosing the optimal value of
$q$ is to use restarting, \cite{nes:07,parnes}. One example of a fixed restart
scheme is as follows: 
\begin{algorithm}[H]
\caption{Fixed restarting}
\label{a-fix_restart}
\begin{algorithmic}[1]
\REQUIRE $x^0 \in \reals^n$, $y^0 = x^0$, $\theta_0= 1$ 
\FOR {$j=0,1,\ldots $}
\STATE carry out Algorithm \ref{a-scheme1} with $q=0$ for $k$ steps
\STATE set $x^0 = x^k$, $y^0 = x^k$ and $\theta_0 = 1$. 
\ENDFOR
\end{algorithmic}
\end{algorithm} 
We restart the algorithm every $k$ iterations, taking
as our starting point the last point produced by the algorithm, where
$k$ is a fixed restart interval. In other words we `forget' all previous
iterations and reset the momentum back to zero.

\paragraph{Optimal fixed restart interval.} 
\label{s-fixres}
We can obtain an upper bound on the optimal restart interval.  If we
restart every $k$ iterations we have, at outer
iteration $j$, inner loop iteration $k$ (just before a restart),
\[ 
f(x^{(j+1,0)}) - f^\star = f(x^{(j,k)}) - f^{\star} \leq 4L\|x^{(j,0)} - x^{\star}\| / k^2 \leq
\left( 8 L  / \mu k^2\right) (f(x^{(j,0)}) - f^\star),
\] 
where the first inequality is the convergence guarantee of Algorithm
\ref{a-scheme1}, and the second comes from the strong convexity of $f$.
So after $jk$ steps we have
\[ 
f(x^{(j,0)}) - f^\star  \leq \left( 8 L / \mu k^2 \right)^j (f(x^{(0,0)}) - f^\star).
\] 
If we assume we have $jk = c$ total iterations and we wish to minimize
$(8L/\mu k^2)^j$ over $j$ and $k$ jointly, we obtain
\begin{equation}
\label{e-fixres}
k^\star =  e \sqrt{8L/\mu}.
\end{equation}
Using this as our restart interval we obtain an accuracy of $\epsilon$ in
less than $\mathcal{O}(\sqrt{L/\mu} \log (1/\epsilon))$
iterations, \ie, the optimal
linear convergence rate as in equation (\ref{e-ordcvg}).

The drawbacks in using fixed restarts are that firstly it depends on unknown
parameters $L$ and, more importantly, $\mu$, and secondly it is a global
parameter that may be inappropriate in better conditioned regions.

%
%
%
%

\subsection{Adaptive restart}
The above analysis suggests that an \emph{adaptive} restart technique may be useful.
In particular we want a scheme that makes some computationally cheap
observation and decides whether or not to restart based on that observation.
In this paper we suggest two schemes 
that perform well in practice and provide some analysis to show accelerated
convergence when these schemes are used.
\BIT \item{\textbf {Function scheme:}} 
we restart whenever 
\[ f(x^k) > f(x^{k-1}). \] 
\item{\textbf {Gradient scheme:}}
we restart whenever
\[
\nabla f(y^{k-1})^T(x^k - x^{k-1}) >0.
\] 
\EIT
Empirically we observe that these two schemes perform similarly well. The
gradient scheme has two advantages over the function scheme. Firstly 
near to the optimum the gradient scheme may be more numerically stable. Secondly
all quantities involved in the gradient scheme are already calculated in accelerated
schemes, so no extra computation is required.

We can give rough justifications for each scheme. The function scheme restarts
at the bottom of the
troughs as in Figure \ref{f-qeff}, thereby avoiding the wasted iterations
where we are moving away from the optimum.  The gradient scheme restarts
whenever the momentum term and the negative gradient
are making an obtuse angle. In other words we restart when the momentum seems to be
taking us in a bad direction, as measured by the negative gradient at that
point. 

\begin{figure}
\begin{center}
\psfrag{k}[t][t]{\small $k$}
\psfrag{fk-fstar}[b][b]{\small $f(x_k) - f^\star$}
\includegraphics[width=0.7\linewidth]{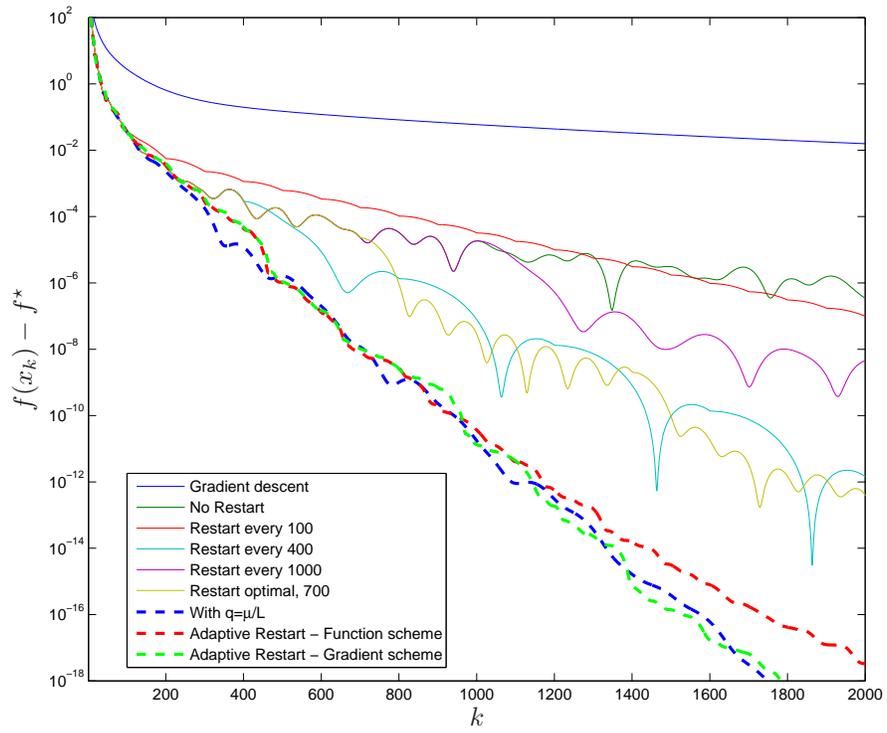}
\caption{Comparison of fixed and adaptive restart intervals.} 
\label{f-adaprestart}
\end{center}
\end{figure}

%
%
Figure \ref{f-adaprestart} shows the effect of different restart intervals
on minimizing a positive definite quadratic function in $n=500$ dimensions. In
this particular case the upper bound on the optimal restart interval is every
$700$ iterations. We note that when this interval is used the convergence is
better than when no restart is used, however not as good as using the optimal
choice of $q$. We also note that restarting every $400$ iterations performs
about as well as restarting every $700$ iterations, suggesting that the optimal
restart interval is somewhat lower than $700$. 
We have also plotted the performance of the two adaptive restart schemes. The
performance is on the same order as the algorithm with the optimal $q$ and much
better than using the fixed restart interval. (Conjugate gradient methods,
\cite{cg}, will generally outperform an accelerated gradient scheme when
minimizing a quadratic; we use quadratics here simply for illustrative purposes.)

\begin{figure} 
\begin{center} 
\psfrag{x1}[t][]{$x_1$}
\psfrag{x2}[b][b]{$x_2$}
\psfrag{q = q-opt}{\footnotesize $q=q^\star$}
\psfrag{q = 0}{\footnotesize $q=0$}
\psfrag{adaptive restart}{\footnotesize adaptive restart}
\includegraphics[width=0.5\linewidth]{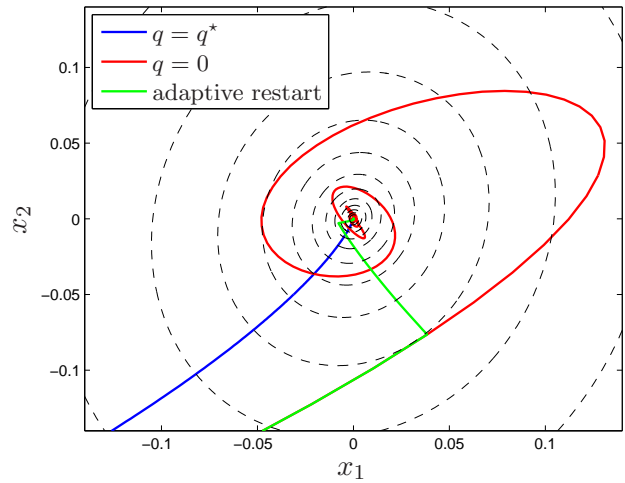} 
\caption{Sequence trajectories under scheme I and with adaptive restart.}
\label{f-trajadap} 
\end{center} 
\end{figure}

Figure \ref{f-trajadap} demonstrates the function restart scheme trajectories
in the two dimensional example, restarting resets the momentum and
prevents the characteristic spiralling behavior.

\section{Analysis}

In this section we consider applying an accelerated scheme to minimizing
a positive definite quadratic.  We shall see that once the momentum is larger
than a critical value we observe periodicity in the iterates. We use this to
prove linear convergence when using adaptive restarting.  The analysis
presented in this section is similar in spirit to the analysis of the heavy
ball method in \cite[\S 3.2]{polyak:87}.

\subsection{Minimizing a quadratic} 
Consider minimizing a convex quadratic. Without
loss of generality we can assume that $f$ has the following form: 
\[ 
f(x) = (1/2)x^TAx 
\] 
where $A \in \reals^{n \times n}$ is positive definite and symmetric. In this
case $x^\star = 0$ and $f^\star = 0$. We have strong convexity
parameter $\mu = \lambda_\mathrm{min} > 0$  
and $L = \lambda_\mathrm{max}$, where $\lambda_\mathrm{min}$ and
$\lambda_\mathrm{max}$ are the minimum and maximum eigenvalues of $A$,
respectively.  

\subsection{The algorithm as a linear dynamical system}
We shall assume a fixed step-size $t=1/L$ for simplicity.
Given quantities $x^{0}$ and $y^{0} = x^{0}$, Algorithm \ref{a-scheme1}
is carried out as follows, 
\[ 
\begin{array}{rcl} x^{k+1} &=& y^{k}
- (1/L) A y^{k} \\ y^{k+1} &=& x^{k+1} + \beta_k (x^{k+1} - x^{k}).
\end{array} 
\]
For the rest of the analysis we shall take $\beta_k$ to be constant and equal to
some $\beta$ for all $k$. This is a somewhat crude approximation, but by making it
we can show that there are two regimes of behavior of the system, depending on
the value of $\beta$.
Consider the eigenvector decomposition of $A = V \Lambda V^T$.  Denote by
$w^{k} = V^T x^{k}$, $v^{k} = V^T y^{k}$. In this basis the update
equations can be written
\[ 
\begin{array}{rcl} w^{k+1} &=& v^{k} - (1/L)
\Lambda v^{k} \\ v^{k+1} &=& w^{k+1} + \beta (w^{k+1} - w^{k}).
\end{array} 
\] 
These are $n$ independently evolving dynamical
systems. The $i$th system evolves according to 
\[ 
\begin{array}{rcl} w_i^{k+1} &=&
v_i^{k} - (\lambda_i/L) v_i^{k} \\ v_i^{k+1} &=& w_i^{k+1} + \beta
(w_i^{k+1} - w_i^{k}), 
\end{array} 
\] 
where $\lambda_i$ is the $i$th eigenvalue of $A$.  Eliminating the sequence
$v_i^{(k)}$ from the above we obtain the following recurrence relation 
for the evolution of $w_i$: 
\[ 
w_i^{k+2} = (1+\beta)(1-\lambda_i/L) w_i^{k+1}
- \beta (1 - \lambda_i /L) w_i^{k},\quad k=0,1,\ldots, 
\] 
where $w_i^{0}$ is known and $w_i^{1} = w_i^{0} (1-\lambda_i /L)$,
$\ie$, a gradient step from $w_i^{0}$. 

The update equation for $v_i$ is identical, differing only in the initial
conditions, 
\[ 
v_i^{k+2} = (1+\beta)(1-\lambda_i /L) v_i^{k+1} - \beta (1 -
\lambda_i /L) v_i^{k},\quad k=0,1,\ldots, 
\] 
where $v_i^{0} = w_i^{0}$ and $v_i^{1} = ((1+\beta)(1-\lambda_i/L)-\beta)v_i^{0}$. 

\subsection{Convergence properties} 
The behavior of this system is determined by the characteristic polynomial of
the recurrence relation,
\begin{equation}
\label{e-poly}
r^2 - (1+\beta)(1-\lambda_i /L) r + \beta (1-\lambda_i /L).
\end{equation}
Let $\beta^\star_i$ be the critical value of $\beta$ for which this polynomial
has repeated roots, \ie,
\[
\beta^\star_i := \frac{1-\sqrt{\lambda_i /L}}{1+\sqrt{\lambda_i /L}}.
\]
If $\beta \leq \beta^\star_i$ then the
polynomial (\ref{e-poly}) has two real roots, $r_1$ and $r_2$, and the system
evolves according to \cite{chiang:84}
\begin{equation}
\label{e-real}
w_i^k = c_1 r_1^k + c_2 r_2^k.  
\end{equation}
When $\beta = \beta^\star_i$ the roots coincide at the point $r^\star = 
(1+\beta)(1-\lambda_i /L)/2 = (1-\sqrt{\lambda_i /L})$; this corresponds to
critical damping. We have the fastest monotone convergence at rate $\propto
(1-\sqrt{\lambda_i /L})^k$. Note that if $\lambda_i = \mu$
then $\beta^\star_i$ is the optimal choice of $\beta$ as given by equation
(\ref{e-bstar}) and the convergence rate is the optimal rate, as given by equation
(\ref{e-lincvg}). This is the case because, as we shall see, the smallest eigenvalue
will come to dominate the convergence of the entire system.

If $\beta < \beta^\star_i$ we are in the low momentum regime, and we say the
system is over-damped. The convergence rate is dominated by the larger root,
which is greater than $r^\star$,\ie, the system exhibits slow monotone convergence.

If $\beta > \beta^\star_i$ then the roots
of the polynomial (\ref{e-poly}) are complex; we are in the high momentum
regime and the system is under-damped and exhibits periodicity. In that case
the characteristic solution is given by \cite{chiang:84}
\[
w_i^{k} = c_i \left(\beta (1-\lambda_i /L)\right)^{k/2} \left(\cos(k
\psi_i - \delta_i)  \right) 
\]
where
\[
\psi_i = \cos^{-1} ((1-\lambda_i /L)(1+\beta)/2\sqrt{\beta (1-\lambda_i /L)}).  
\] 
and $\delta_i$ and $c_i$ are constants that depend on the initial conditions;
in particular for $\beta \approx 1$ we have $\delta_i \approx 0$ and
we will ignore it. Similarly,
\[ 
v_i^{k} = \hat c_i \left(\beta (1-\lambda_i /L)\right)^{k/2} \left(\cos(k
\psi_i - \hat \delta_i)  \right) 
\] 
where $\hat \delta_i$ and $\hat c_i$ are constants, and again $\hat \delta_i
\approx 0$.
For small
$\theta$ we know that $\cos^{-1}(\sqrt{1-\theta}) \approx \sqrt{\theta}$, and
therefore if $\lambda_i \ll L$, then 
\[ 
\psi_i \approx \sqrt{\lambda_i /L}.  
\] 
In particular the frequency of oscillation for the mode corresponding to the
smallest eigenvalue $\mu$ is approximately given by $\psi_\mu \approx \sqrt{\mu/L}$. 

To summarize, based on the value of $\beta$ we observe the following behaviors: 
\BIT
\setlength{\itemsep}{-3pt}
\item $\beta > \beta_i^\star$: high momentum, under-damped
\item $\beta < \beta_i^\star$: low momentum, over-damped
\item $\beta = \beta_i^\star$: optimal momentum, critically damped.
\EIT

\subsection{Observable quantities} 
We don't observe the evolution of the modes, but we can observe the evolution
of the function value; which is given by
\[
f(x^{k}) = \sum_{i=1}^n  (w_i^{k})^2 \lambda_i. 
\] 
and if $\beta > \beta^\star = (1-\sqrt{\mu /L})/(1+\sqrt{\mu /L})$
we are in the high momentum regime for all modes and thus
\[ 
f(x^{k}) = \sum_{i=1}^n  (w_i^{k})^2 \lambda_i \approx \sum_{i=1}^n
(w_i^{0})^2\lambda_i \beta^k (1-\lambda_i /L)^k \cos^2(k \psi_i).  
\] 
The function value will quickly be dominated by the smallest eigenvalue
and we have that 
\begin{equation}
\label{e-fctsch}
f(w^{k}) \approx (w_\mu^{0})^2 \mu \beta^k (1-\mu /L)^k \cos^2\left(k \sqrt{\mu/L}\right), 
\end{equation}
where we have replaced $\psi_\mu$ with $\sqrt{\mu /L}$, and we are using the
subscript $\mu$ to denote those quantities corresponding to that mode.
A similar analysis for the gradient restart scheme yields
\begin{equation}
\label{e-grasch}
\nabla f(y^k)^T(x^{k+1} - x^k) \approx \mu v^k_\mu(w^{k+1}_\mu - w^k_\mu) \propto \beta^k(1-\mu/L)^k \sin(2k\sqrt{\mu/L}).
\end{equation}
In other words observing the quantities in  (\ref{e-fctsch}) or (\ref{e-grasch})
we expect to see oscillations at a frequency proportional to
$\sqrt{\mu/L}$, \ie, the frequency of oscillation is telling
us something about the condition number of the function. 

\subsection{Convergence with adaptive restart} 
If we apply Algorithm \ref{a-scheme1} with $q=0$ 
to minimize a quadratic we start with $\beta_0 = 0$, \ie, the system is in the
low momentum, monotonic regime. Eventually $\beta_k$ becomes larger than
$\beta^\star$ and we enter the high momentum, oscillatory regime. It takes about
$(3/2) \sqrt{L/\mu}$ iterations for $\beta_k$ to exceed $\beta^\star$. After that
the system is under-damped and the iterates obey equations (\ref{e-fctsch}) and 
(\ref{e-grasch}).  Under either adaptive restart scheme, equations (\ref{e-fctsch})
and (\ref{e-grasch}) indicate that we shall observe the restart condition
after a further $(\pi/2)\sqrt{L/\mu}$ iterations.  We restart and the process
begins again, with $\beta_k$ set back to zero. Thus under either scheme we
restart approximately every
\[
k^\star = \frac{\pi+3}{2} \sqrt{\frac{L}{\mu}}
\]
iterations (\cf, the upper bound on optimal fixed restart interval
(\ref{e-fixres})).  Following a similar derivation to \S \ref{s-fixres}, 
this restart interval guarantees us an accuracy of $\epsilon$ within 
$\mathcal{O}(\sqrt{L/\mu} \log (1/\epsilon))$
iterations,
\ie, we have recovered the optimal linear convergence rate of equation
(\ref{e-ordcvg}) via adaptive restarting, with no prior knowledge of $\mu$.

\subsection{Extension to smooth convex minimization} 
In many cases the function we are minimizing is well approximated by a
quadratic near the optimum, \ie, there is a region inside of which
\[
f(x) \approx f(x^\star) + (x- x^\star)^T \nabla^2 f(x^\star) (x-x^\star),
\]
and loosely speaking we are minimizing a quadratic. Once we
are inside this region we will observe behavior consistent with the
analysis above, and we can exploit this behavior to achieve fast convergence by
using restarts.  Note that the Hessian at the optimum may have smallest
eigenvalue $\lambda_\mathrm{min} > \mu$, the global strong convexity parameter,
in other words we can achieve a faster local convergence than even if we had
exact knowledge of the global parameter.
This result is similar in spirit to the restart method applied to the non-linear
conjugate gradient method, where it is desirable to restart the algorithm once
it reaches a region in which the function is well approximated by a quadratic
\cite[\S 5.2]{nocedal}. 

The effect of these restart schemes outside of the quadratic region is unclear.
In practice we observe that restarting based on one of the criteria described
above is almost always helpful, even far away from the
optimum. However, we have observed cases where restarting far from the optimum
can slow down the early convergence slightly, until the quadratic region is
reached and the algorithm enters the rapid linear convergence phase. 



\section{Numerical examples}
In this section we describe three further numerical examples that demonstrate
the improvement of accelerated algorithms under an adaptive restarting technique.

\subsection{Log-sum-exp}
Here we minimize a smooth convex function that is \emph{not} strongly convex.
Consider the following optimization problem
\[
\begin{array}{ll}
\mbox{minimize} & \rho\log\left(\sum_{i=1}^m \exp\left((a_i^T x - b_i)/\rho\right) \right)
\end{array}
\]
where $x \in \reals^n$.
The objective function is smooth, but not strongly convex, it
grows linearly asymptotically. Thus, the optimal value of $q$ in Algorithm
\ref{a-scheme1} is zero.  The quantity $\rho$ controls the smoothness of the
function, as $\rho \rightarrow 0$, $f(x) \rightarrow \max_{i=1,\ldots,m}
(a_i^Tx-b_i)$.  As it is smooth, we expect the 
region around the optimum to be well approximated by a quadratic (we
consider only examples where the optimal value is finite), and thus
we expect to eventually enter a region where our restart method will
obtain linear convergence without any knowledge of where this region is,
the size of the region or the local function parameters within this region.
For smaller values of $\rho$ the smoothness of the objective function decreases
and thus we expect to take more iterations before we
enter the region of linear convergence.

\begin{figure}[h]
\begin{center}
\psfrag{k}[t][t]{\small $k$}
\psfrag{fk-fstar}[b][b]{\small $(f_k - f^\star)/f^\star$}
\psfrag{rho=0.05}{$\rho=0.05$}
\psfrag{rho=0.10}{$\rho=0.1$}
\psfrag{rho=0.50}{$\rho=0.5$}
\psfrag{rho=1.00}{$\rho=1$}
\includegraphics[width=1\linewidth]{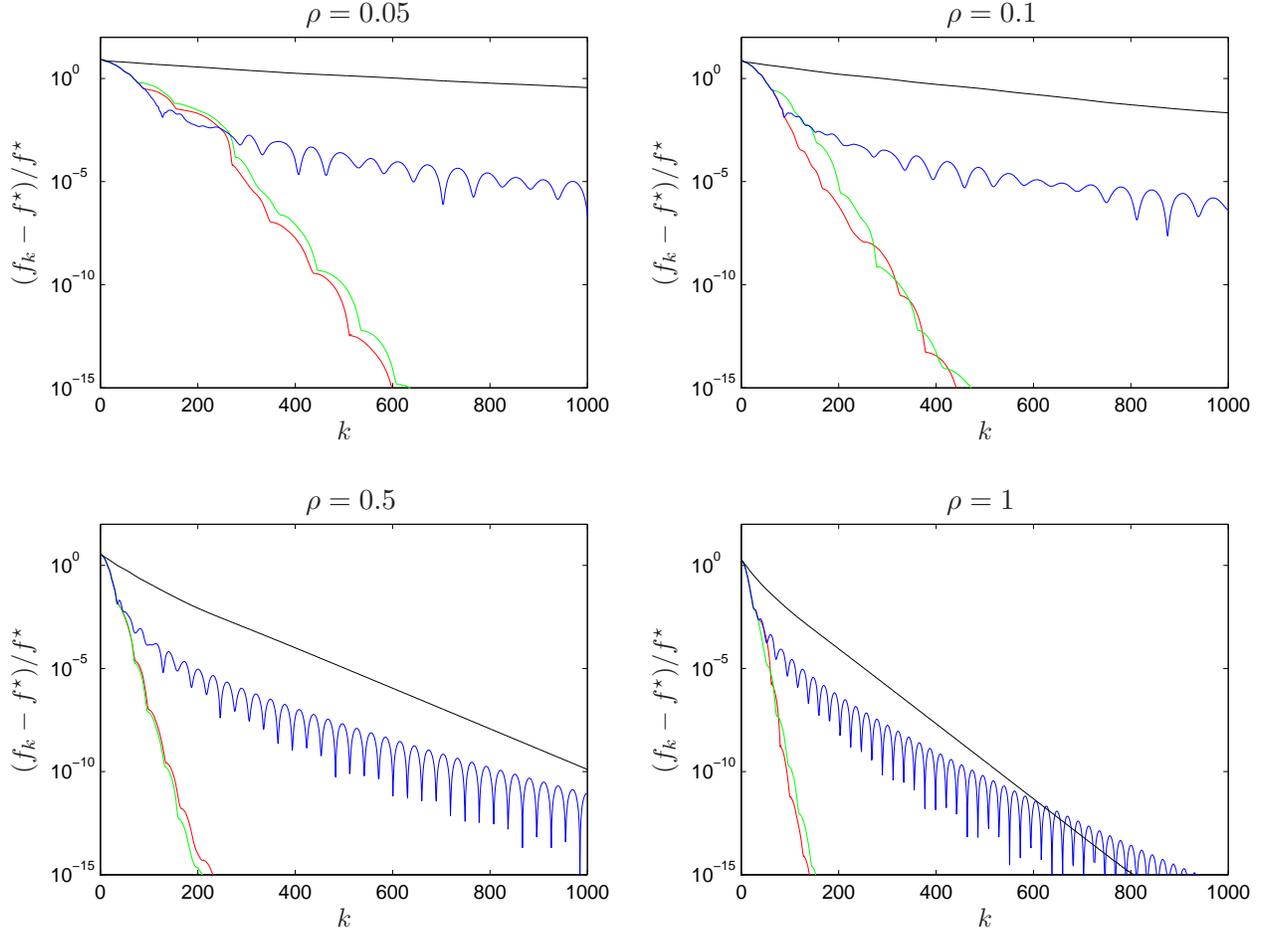}
\caption{Minimizing a smooth but not strongly convex function; the black line
is gradient descent, the blue line is Algorithm \ref{a-scheme1}, the red line
is the function adaptive restart scheme, the green line is the gradient
adaptive restart scheme.} 
\label{f-logsumexp}
\end{center}
\end{figure}

As a particular example we took $n=20$ and $m=100$; we generated the
$a_i$ and $b_i$ randomly. Figure \ref{f-logsumexp} demonstrates the
performance of four different schemes for four different values of
$\rho$. We selected the step size for each case using backtracking. We
note that both restart schemes perform well, eventually beating both
gradient descent and the accelerated scheme.  Both the function and
gradient schemes eventually enter a region of fast linear convergence.
For large $\rho$ we see that even gradient descent performs well, as,
similar to the restarted method, it is able to automatically exploit
the local strong convexity of the quadratic region around the optimum.
Notice also the appearance of the periodic behavior.

\subsection{Sparse linear regression}
Consider the following optimization problem:
\begin{equation}
\label{e-lasso}
\begin{array}{ll}
\mbox{minimize} & (1/2)\|Ax-b\|_2^2 + \rho \|x\|_1,
\end{array}
\end{equation}
over $x\in \reals^n$, where $A \in \reals^{m\times n}$ and in general
$n \gg m$.
This is a widely studied problem in the field of compressed sensing, see \eg,
\cite{candes:06,donoho:06,candes:08,tibshirani:94}. Loosely speaking problem
(\ref{e-lasso}) seeks a sparse vector with a small measurement error. 
The quantity $\rho$ trades off these two competing objectives. The iterative
soft-threshold algorithm (ISTA) can be used to solve (\ref{e-lasso})
\cite{chambolle:98,daub:04}. ISTA relies on
the soft-thresholding operator:
\[
\mathcal{T}_\alpha (x) = \mathrm{sign}(x) \max(|x|-\alpha,0),
\]
where all the operations are applied elementwise.
The ISTA algorithm, with constant step-size $t$, is given by
\begin{algorithm}[H]
\caption{ISTA}
\label{a-ista}
\begin{algorithmic}[1]
\REQUIRE $x^{(0)} \in \reals^n$
\FOR {$k=0,1,\ldots$}
\STATE $x^{k+1} = \mathcal{T}_{\rho t} (x^k - tA^T(Ax^k - b))$.
\ENDFOR
\end{algorithmic}
\end{algorithm} 
The convergence rate of ISTA is guaranteed to be at least $\mathcal O(1/k)$,
making it analogous to gradient descent.

The fast iterative soft thresholding algorithm (FISTA) was developed in
\cite{fista}; a similar algorithm was also developed by Nesterov in
\cite{nes:07}. FISTA essentially applies acceleration to the ISTA algorithm; it
is carried out as follows, 
\begin{algorithm}[H]
\caption{FISTA}
\label{a-fista}
\begin{algorithmic}[1]
\REQUIRE $x^{(0)} \in \reals^n$, $y^0 = x^0$ and $\theta_0 = 1$
\FOR {$k=0,1,\ldots$}
\STATE $x^{k+1} = \mathcal{T}_{\rho t} (y^k - tA^T(Ay^k - b))$
\STATE $\theta_{k+1} = (1+\sqrt{1+4\theta_{k}^2})/2$ 
\STATE $\beta_{k+1}=  (\theta_k-1)/\theta_{k+1}$
\STATE $y^{k+1} = x^{k+1} + \beta_{k+1}(x^{k+1} - x^k)$.
\ENDFOR
\end{algorithmic}
\end{algorithm}
For any choice of $t \leq 1/\lambda_\mathrm{max}(A^T A)$
FISTA obtains a convergence rate of at least $O(1/k^2)$.
The objective in problem (\ref{e-lasso}) is non-smooth, so it does not fit the
class of problems we are considering in this paper. However we are
seeking a sparse solution vector $x$, and we note that once the
non-zero basis of the solution has been identified we are essentially 
minimizing a quadratic.
Thus we expect that after a certain number of iterations
adaptive restarting can provide linear convergence.

It is easy to show that the function adaptive restart scheme can be performed
without an extra application of the matrix $A$, which is the costly
operation in the algorithm.

In performing FISTA we do not evaluate a gradient, however FISTA can be thought
of as a \emph{generalized gradient} scheme, in which we take
\[
x^{k+1} = \mathcal{T}_{\lambda t} (y^k - tA^T(Ay^k - b)) := y^k - tG(y^k)
\]
to be a generalized gradient step, where $G(y^k)$ is the generalized gradient
at $y^k$.  In this case the gradient restart scheme amounts to restarting
whenever
\[
G(y^k)^T(x^{k+1}-x^k) >0,
\]
or equivalently
\begin{equation}
\label{e-gengrad}
(y^k - x^{k+1})^T(x^{k+1} - x^k)>0.
\end{equation}
%
%

\begin{figure}
\begin{center}
\psfrag{k}[t][t]{\small $k$}
\psfrag{fk-fstar}[b][b]{\small $(f(x^k) - f^\star)/f^\star$}
\psfrag{m100}[b][b]{\small $n=2000$, $m=100$, $s=20$}
\psfrag{m500}[b][b]{\small $n=2000$, $m=500$, $s=100$}
\includegraphics[width=1\linewidth]{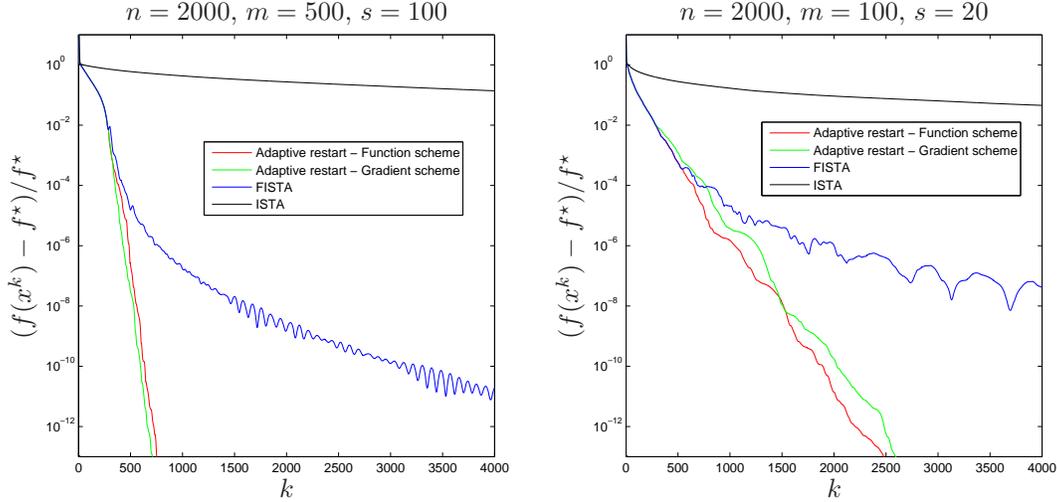}
\caption{Adaptive restarting applied to the FISTA algorithm.}
\label{f-lasso_one_fig}
\end{center}
\end{figure}

We generated data for the numerical instances as follows. Firstly the entries
of $A$ were sampled from a standard normal distribution.  We then randomly
generated a sparse vector $y$ with $n$ entries, only $s$ of which were
non-zero.  We then set $b = Ay + w$, where the entries in $w$ were IID sampled
from $\mathcal{N}(0,0.1)$.  This ensured that the solution vector $x^\star$ is
approximately $s$-sparse. We chose $\rho =1$ and the step size
$t = 1/\lambda_\mathrm{max}(A^TA)$.
Figure \ref{f-lasso_one_fig} shows the dramatic speedup that adaptive
restarting can provide, for two different examples.  

\subsection{Quadratic programming}
Consider the following quadratic program,
\begin{equation}
\label{e-qp}
\begin{array}{ll}
\mbox{minimize} & (1/2) x^T Q x + q^T x \\
\mbox{subject to} & a \leq x \leq b,
\end{array}
\end{equation}
over $x \in \reals^n$, where $Q \in \reals^{n \times n}$ is positive definite
and $a,b \in \reals^n$ are fixed vectors. The constraint inequalities
are to be interpreted element-wise, and we assume that $a < b$.
We denote by $\Pi_\mathcal{C}(z)$ the projection of a point $z$ onto the
constraint set, which amounts to thresholding the entries in $z$.

Projected gradient descent can solve (\ref{e-qp}); it is carried out as
follows,
\[
x^{k+1} = \Pi_\mathcal{C} \left(x^{k} - t (Qx^k +q) \right).
\]
Projected gradient descent obtains a guaranteed convergence rate of $\mathcal{O}(1/k)$.
Acceleration has been successfully applied to the projected gradient method,
\cite{nes:07,fista}.
\begin{algorithm}[H]
\caption{Accelerated projected gradient}
\label{a-qp}
\begin{algorithmic}[1]
\REQUIRE $x^{0} \in \reals^n$, $y^{0} = x^{0}$ and $\theta_0=1$
\FOR {$k=0,1,\ldots $}
\STATE $x^{k+1} = \Pi_\mathcal{C} \left(y^{k} - t (Qy^k +q) \right)$
\STATE $\theta_{k+1}$ solves $\theta_{k+1}^2 =
(1-\theta_{k+1})\theta_k^2$
\STATE $\beta_{k+1} = \theta_k(1 -
\theta_k)/(\theta_k^2 + \theta_{k+1})$ 
\STATE $y^{k+1} = x^{k+1} + \beta_{k+1}(x^{k+1} - x^{k})$
\ENDFOR
\end{algorithmic}
\end{algorithm} 
For any choice of $t \leq 1/\lambda_\mathrm{max}(Q)$ accelerated projected
gradient schemes obtain a convergence rate of at least $\mathcal{O}(1/k^2)$.

The presence of constraints make this a non-smooth optimization problem,
however once the constraints that are active have been identified 
the problem reduces to minimizing a quadratic on a subset of the variables, and
we expect adaptive restarting to increase the rate of convergence. As in the
sparse regression example we can use the
generalized gradient in our gradient based restart scheme, \ie, we restart
based on condition (\ref{e-gengrad}).

\begin{figure}[h]
\begin{center}
\psfrag{k}[t][t]{\small $k$}
\psfrag{fk-fstar}[b][b]{\small $(f(x^k) - f^\star)/f^\star$}
\includegraphics[width=0.7\linewidth]{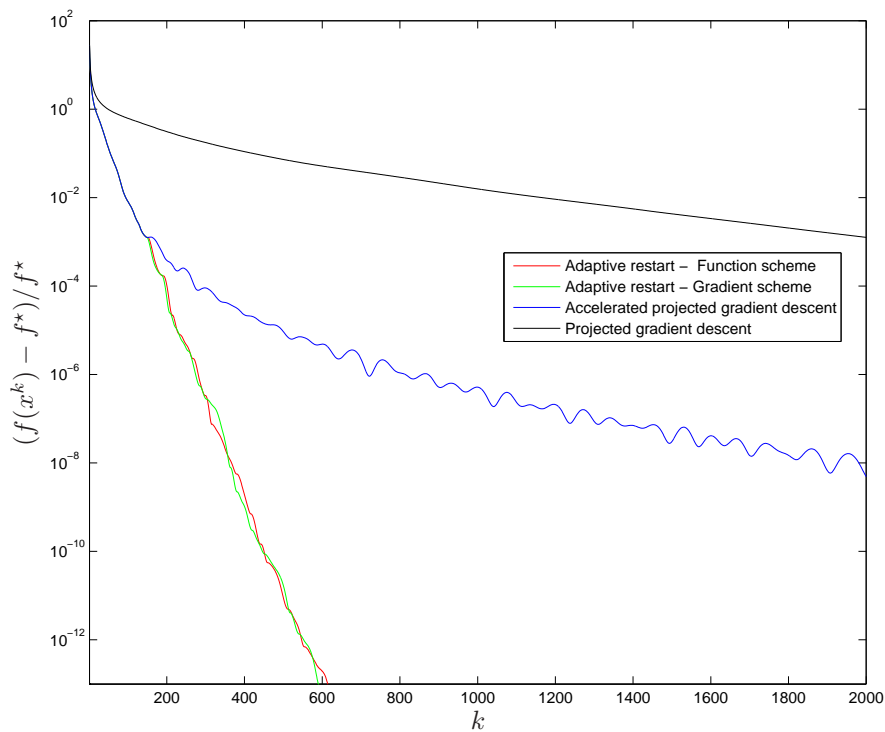}
\caption{Adaptive restarting applied to the accelerated projected gradient algorithm.}
\label{f-qp}
\end{center}
\end{figure}

As a final example, we set $n=500$ and generate $Q$ and $q$ randomly;
$Q$ has a condition number of $10^7$. We take $b$ to be the vector of
all ones, and $a$ to be that of all negative ones. For information,
the solution to this problem has $70$ active constraints. The
step-size is set to $t=1/\lambda_\mathrm{max}(Q)$ for all algorithms.
Figure \ref{f-qp} shows the performance of projected gradient descent,
accelerated projected gradient descent, and the two restart
techniques.

\section{Summary}
In this paper we have demonstrated a simple heuristic adaptive restart
technique that can improve the convergence performance of accelerated
gradient schemes for smooth convex optimization.  We restart the algorithm
whenever we observe a certain condition on the objective function value or
gradient value.  We provided some qualitative analysis to show that we can
recover the optimal linear rate of convergence in many cases; in particular
near the optimum of a smooth function we can potentially dramatically
accelerate the rate of convergence, even if the function is not globally
strongly convex.  We demonstrated the performance of the scheme on some simple
numerical examples.

\section*{Acknowledgments}
We are very grateful to Stephen Boyd for his help and
encouragement. We would also like to thank Stephen Wright for his
advice and feedback, and Stephen Becker and Michael Grant for useful
discussions. E.~C.~would like to thank the ONR (grant
N00014-09-1-0258) and the Broadcom Foundation for their support.

\bibliographystyle{plain} 
\bibliography{adap_restart_paper}
\end{document}

%% file: defs.tex
\newcommand{\BEAS}{\begin{eqnarray*}}
\newcommand{\EEAS}{\end{eqnarray*}}
\newcommand{\BEA}{\begin{eqnarray}}
\newcommand{\EEA}{\end{eqnarray}}
\newcommand{\BEQ}{\begin{equation}}
\newcommand{\EEQ}{\end{equation}}
\newcommand{\BIT}{\begin{itemize}}
\newcommand{\EIT}{\end{itemize}}
\newcommand{\BNUM}{\begin{enumerate}}
\newcommand{\ENUM}{\end{enumerate}}

\newcommand{\BA}{\begin{array}}
\newcommand{\EA}{\end{array}}

\newcommand{\cf}{{\it cf.}}
\newcommand{\eg}{{\it e.g.}}
\newcommand{\ie}{{\it i.e.}}

\newcommand{\reals}{{\mbox{\bf R}}}

{\begin{quote}}{\end{quote}}

\makeatletter
\long\def\@makecaption#1#2{
   \vskip 9pt 
   \begin{small}
   \setbox\@tempboxa\hbox{{\bf #1:} #2}
   \ifdim \wd\@tempboxa > 5.5in
        \begin{center}
        \begin{minipage}[t]{5.5in}
        \addtolength{\baselineskip}{-0.95pt}
        {\bf #1:} #2 \par
        \addtolength{\baselineskip}{0.95pt}
        \end{minipage}
        \end{center}
   \else 
	\hbox to\hsize{\hfil\box\@tempboxa\hfil}  
   \fi
   \end{small}\par
}
\makeatother

\newcounter{oursection}

\newcounter{lecture}